\documentclass{amsart}

\usepackage[dvips]{graphicx}
\usepackage{amscd}
\usepackage{amsmath}
\usepackage{amsfonts}
\usepackage{amssymb}

\usepackage[all]{xy}

\def\openone
{\mathchoice
{\hbox{\upshape \small1\kern-3.3pt\normalsize1}}
{\hbox{\upshape \small1\kern-3.3pt\normalsize1}}
{\hbox{\upshape \tiny1\kern-2.3pt\SMALL1}}
{\hbox{\upshape \Tiny1\kern-2pt\tiny1}}}

\newtheorem{theorem}{Theorem}
\newtheorem{corollary}[theorem]{Corollary}
\newtheorem{lemma}[theorem]{Lemma}
\newtheorem{proposition}[theorem]{Proposition}
\theoremstyle{definition}

\theoremstyle{remark}

\def\1{{\openone}}
\def\A{{\mathfrak{A}}}

\def\N{{\bf N}}

\def\T{{\bf T}}
\def\Z{{\bf Z}}
\def\C{{\bf C}}

\def\Ad{\mathop{\rm Ad}\nolimits}

\def\Aut{\mathop{\rm Aut}\nolimits}
\def\UHF{{\rm UHF}}
\def\supp{\mathop{\rm supp}\nolimits}
\def\End{\mathop{\rm End}\nolimits}

\title[Homogeneity of the pure state space of the Cuntz algebra]{Homogeneity
of the pure state space of the Cuntz algebra}
\author{Ola~Bratteli}
\address{Mathematics Institute\\
University of Oslo\\
PB 1053 Blindern\\
N-0316 Oslo\\ Norway}
\author{Akitaka~Kishimoto}
\address{Department of Mathematics\\
Hokkaido University\\
Sapporo\\ 060 Japan}

\begin{document}

\maketitle

\begin{abstract}
If $\omega_1,\omega_2$ are two pure gauge-invariant states of the Cuntz
algebra
${\mathcal O}_d$, we show that there is an automorphism $\alpha$ of
${\mathcal O}_d$ such that $\omega_1=\omega_2\circ\alpha$. If $\omega$ is a
general pure state on ${\mathcal O}_d$ and $\varphi_0$ is a given Cuntz
state, we show that there exists an endomorphism $\alpha$ of ${\mathcal O}_d$
such that $\varphi_0=\omega\circ\alpha$
\end{abstract}

\section{Introduction}

Let $\A$ be a simple separable C*-algebra, and let $\pi_1, \pi_2$ be
representations of $\A$ on Hilbert spaces ${\mathcal H}_1,{\mathcal H}_2$. The
representations
$\pi_1,\pi_2$ are said to be algebraically equivalent if $\pi_1(\A)''$ and
$\pi_2(\A)''$ are isomorphic von Neumann algebras. If there is an
automorphism $\alpha$ of $\A$ such that $\pi_1$ and $\pi_2\circ \alpha$ are
quasi-equivalent, then $\pi_1,\pi_2$ are clearly algebraically equivalent.
Powers proved in \cite{Pow67} that if $\A$ is a UHF algebra the converse is
true. His method extends readily to the case that $\A$ is an AF-algebra,
\cite{Bra72}.
See also section~12.3 in \cite{KR86}. In the special case that $\pi_1$  (and
therefore $\pi_2)$ is irreducible, Kadison's transitivity theorem
therefore implies that if $\A$ is a simple AF algebra and if
$\omega_1$ and $\omega_2$ are pure states on $\A$, there exists an
automorphism $\alpha$ of $\A$ such that $\omega_1=\omega_2\circ
\alpha$. To our knowledge, this question has only been settled in the
affirmative when
$\A$ is an AF-algebra. As a beginning of a possible resolution of the question
for purely infinite algebras, we here prove the statements in the abstract.
Recall from \cite{Cun77} that the Cuntz algebra
${\mathcal O}_d$ is the C*-algebra generated by $d$ operators $s_1,\ldots,s_d$
satisfying
\begin{eqnarray*}
&&s_j^\ast s_i=\delta_{ij}\1 \\
&&\sum_{i=1}^d s_is_i^\ast=\1
\end{eqnarray*}
There is an action $\gamma$ of the group $U(d)$ of unitary $d\times d$
matrices on
${\mathcal O}_d$ given by
$$
\gamma_g(s_i)=\sum_{j=1}^d g_{ji}s_j
$$
for $g=[g_{ij}]_{i,j=1}^d$ in $U(d)$. In particular the {\em gauge action\/}
$\tau=\gamma|_{\T}$ is defined by
$$
\tau_z(s_i)=zs_i\;,\qquad z\in\T\subset\C\;.
$$
If $\UHF_d$ is the fixed point subalgebra under the gauge action, then
$\UHF_d$ is the closure of the linear span of all Wick ordered polynomials of
the form
$$
s_{i_1}\ldots s_{i_k} s_{j_k}^\ast\ldots s_{j_1}^\ast
$$
$\UHF_d$ is isomorphic to the UHF algebra of Glimm type $d^\infty$:
$$
\UHF_d\cong M_{d^\infty}=\bigotimes_1^\infty M_d
$$
in such a way that the isomorphism carries the Wick ordered polynomial above
into the matrix element
$$
e_{i_1 j_1}^{(1)}\otimes e_{i_2 j_2}^{(2)}\otimes\cdots\otimes
e_{i_k j_k}^{(k)}\otimes\1\otimes\1\otimes\cdots\;.
$$
In the case that $d$ is a power of a prime, the gauge action $\tau$ is 
in fact characterized by the fact that its fixed
point algebra is isomorphic to $\UHF_d$, i.e. if $\alpha$ is another faithful
action of $\T$  on ${\mathcal O}_d$ such that the fixed point algebra
${\mathcal O}_d^\alpha$ is isomorphic to $\UHF_d$, then either
$z\mapsto\alpha_z$ or $z\mapsto\alpha_z^{-1}$ is conjugate to $\tau$. This
follows from \cite[Corollary~4.1]{BK99}. (Since $\UHF_d$ is simple and
$\alpha$ is faithful, the crossed product ${\mathcal O}_d\times_\alpha \T$ is
stably isomorphic to $\UHF_d$, \cite{KT78}, and in particular it is simple.
Since
$$
{\mathcal O}_d^\alpha\cong P_\alpha(0)({\mathcal O}_d
\times_\alpha\T)P_\alpha(0)\;,
$$

\noindent
$[P_\alpha(0)]$ is just $[\1]$ when $K_0({\mathcal O}_d\times_\alpha\T)$ is
identified with $K_0({\mathcal O}_d^\alpha)$. By the Pimsner-Voiculescu exact
sequence it follows that
$\widehat{\alpha}_\ast$ on
$K_0({\mathcal O}_d\times_\alpha\T)=\Z[\frac{1}{d}]$ is multiplication by $d$
or $1/d$. For this last argument it is important that $d$ is a power of a
prime, as seen from the example $d=6$ and $\widehat{\alpha}_\ast$ 
equal to multiplication by 4/9 on $\Z[\frac{1}{6}]$) Because of this, our main result Theorem~5 can be given the
following more universal form:

\begin{corollary}
Assume that $d$ is a power of a prime. Let $\varphi_1$ and 
$\varphi_2$ be pure states on ${\mathcal O}_d$, and
assume that there exist actions $\alpha_i$ of $\T$ on ${\mathcal O}_d$ 
such that
${\mathcal O}_d^{\alpha_i}\cong\UHF_d$ and
$\varphi_i\circ\alpha_i=\varphi_i$ for $i=1,2$. Then there exists an
automorphism $\beta$ of ${\mathcal O}_d$ such that
$$
\varphi_1=\varphi_2\circ\beta
$$
\end{corollary}

The question whether any pure state on ${\mathcal O}_d$ is invariant under a
gauge action is left open.

The restriction of $\gamma_g$ to $\UHF_d$
is carried into the action
$$
\Ad(g)\otimes\Ad(g)\otimes\cdots
$$
on $\bigotimes\limits_1^\infty M_d$. We define the canonical endomorphism
$\lambda$ on $\UHF_d$ (or on ${\mathcal O}_d$) by
$$
\lambda(x)=\sum_{j=1}^d s_j x s_j^\ast
$$
and the isomorphism carries $\lambda$ over into the one-sided shift
$$
x_1\otimes x_2\otimes x_3\otimes \cdots \to \1\otimes x_1\otimes
x_2\otimes\cdots
$$
on $\bigotimes\limits_1^\infty M_d$.

If $\eta_1,\ldots,\eta_d$ are complex scalars with
$\sum\limits_{j=1}^d|\eta_j|^2=1$, we can define a state on ${\mathcal O}_d$
by
$$
\varphi_\eta(s_{i_1}\ldots s_{i_k}\ s_{j_\ell}^\ast\ldots s_{j_1}^\ast)
=\eta_{i_1}\ldots \eta_{i_k}\ \overline{\eta_{j_\ell}}\ldots
\overline{\eta_{j_1}}
$$
\cite{Cun77}, \cite{Eva80}, \cite{BJP96}, \cite{BJ97},
\cite{BJKW}.

This state is pure, and non-gauge invariant, and the $U(d)$ action is
transitive on these states, which are called Cuntz states. The restriction of
$\varphi_\eta$ to $\UHF_d$ identifies with the pure product state given by
infinitely many copies of the vector state defined by the vector
$(\eta_1,\ldots,\eta_d)$ on
$M_d$.

In this paper we will also consider the one-one correspondence between the
set
${\mathcal U}({\mathcal O}_d)$ of unitaries in ${\mathcal O}_d$ and the set
${\rm End}({\mathcal O}_d)$ of unital endomorphisms of ${\mathcal O}_d$. If
$u\in{\mathcal U}({\mathcal O}_d)$ then $\alpha_u(s_i)=us_i$ defines an
endomorphism, and if $\alpha\in{\rm End}({\mathcal O}_d)$ the corresponding
unitary is $u=\sum\limits_{i=1}^d \alpha(s_i)s_i^\ast$. It has been proved by
R\o{}rdam that
$$
{\mathcal U}_i=\{u\in{\mathcal U}({\mathcal O}_d)|\alpha_u\
\mbox{ is an inner automorphism}\}
$$
is a dense subset of ${\mathcal U}({\mathcal O}_d)$,
[R\o{}r93].
We give a
shorter proof of this, and also show that
$$
{\mathcal U}_a=\{u\in{\mathcal U}({\mathcal O}_d)|\alpha_u\
\mbox{ is an automorphism}\}
$$
is a dense $G_\delta$ subset of ${\mathcal U}({\mathcal O}_d)$ such that the
complement ${\mathcal U}({\mathcal O}_d)\setminus {\mathcal U}_a$ is also
dense.

By using the above correspondence between ${\mathcal U}({\mathcal O}_d)$ and
$\End({\mathcal O}_d)$, it follows (see the proof of Proposition~8) that if
$\omega$ is a pure state and $\varphi_0$ a Cuntz state there exists an
endomorphism $\alpha$ of ${\mathcal O}_d$ such that
$\varphi_0=\omega\circ\alpha$. Although the automorphism group is dense in
$\End({\mathcal O}_d)$ (in the topology of pointwise convergence), the
question whether $\alpha$ can be chosen to be an automorphism is left open (in
this approach).

\section{Transitivity of the automorphism group on the pure gauge-invariant
states}

\setcounter{equation}{0}

In this section we prove the first main result mentioned in the abstract.

Let $\UHF_d$ be the UHF algebra of type $d^\infty$ and let $(A_n)$ be an
increasing sequence of C*-subalgebras of $\UHF_d$ such that
$\UHF_d=\overline{\cup A_n}$ and $A_n\cong M_{d^n}$. We first use Power's
transitivity on $\UHF_d$ to find an approximate factorization for any pure
state on $\UHF_d$:

\begin{lemma}
Let $\varphi$ be a pure state of $\UHF_d$ and $\varepsilon>0$. Then there
exists a pure state $\varphi'$ of $\UHF_d$, an increasing sequence $\{B_n\}$
of finite type I subfactors of $\UHF_d$, and an increasing subsequence
$\{k_n\}$ in
$\N$ such that $\varphi'|B_n$ is a pure state of $B_n$ and $A_{k_n}\subset
B_n\subset A_{k_{n+1}}$ for every $n$, and
$$
\|\varphi-\varphi'\|<\varepsilon\;.
$$
\end{lemma}

\begin{proof}
Since the automorphism group $\Aut(\UHF_d)$ of $\UHF_d$ acts transitively
on the set of pure states of $\UHF_d$, \cite{Pow67}, there exists an
increasing sequence
$\{D_n\}$ of finite type I subfactors of $\UHF_d$ such that $D_n\cong
M_{d^n}$ and
$\varphi|D_n$ is pure for every $n$. Then we can find sequences $\{u_n\}$ and
$\{v_n\}$ of unitaries in $\UHF_d$ and increasing sequences $\{k_n\}$ and
$\{\ell_n\}$ in $\N$ such that
\begin{eqnarray*}
&&A_{k_1}\subset\Ad(v_1u_1)(D_{\ell_1})\subset A_{k_2}\subset
         \Ad(v_2u_2v_1u_1)(D_{\ell_2})\subset A_{k_3}\subset \cdots \\
&&u_n\in \UHF_d\cap\Ad(v_{n-1}u_{n-1}\ldots v_1u_1)(D_{\ell_{n-1}})' \\
&&v_n\in \UHF_d\cap A_{k_n}' \\
&&\|u_n-1\|< \varepsilon/2^{n+2}\qquad
\|v_n-1\|<\varepsilon/2^{n+2}
\end{eqnarray*}
where $D_0=\C1$. (Let $k_1=1$. Then we choose $u_1$ and $\ell_1$ such that
$A_{k_1}\subset\Ad u_1(D_{\ell_1})$ and $\|u_1-1\|<\varepsilon/8$. Further we
choose $k_2$ and $v_1$ such that $v_1\in \UHF_d\cap A_{k_1}'$,
$\|v_1-1\|<\varepsilon/8$, and, $\Ad(v_1u_1)(D_{\ell_1})\subset A_{k_2}$. We
just repeat this process.) Then the limit $w=\lim v_n u_n\ldots v_1u_1$ exists
and is a unitary such that $\|w-1\|<\varepsilon/2$ and
$$
A_{k_1}\subset\Ad w(D_{\ell_1})\subset A_{k_2}\subset \Ad
w(D_{\ell_2})\subset\cdots
$$
Let $\varphi'=\varphi\circ\Ad w^\ast$. Then $\varphi'$ is a pure state with
$\|\varphi-\varphi'\|<\varepsilon$ and $\varphi'|\Ad w(D_{\ell_n})$ is a pure
state for every $n$. Put $B_n=\Ad w(D_{\ell_n})$.
\end{proof}

We next show that for any pair of pure states $\varphi_1,\varphi_2$ on
$\UHF_d$, there is a tensor product decomposition of $\UHF_d$ such that
$\varphi_1,\varphi_2$ have approximate factorizations with respect to certain
sub-decompositions (necessarily different for $\varphi_1$ and $\varphi_2$):

\begin{lemma}
Let $\varphi_1$ and $\varphi_2$ be pure states of $\UHF_d$ and let
$\varepsilon>0$. Then there exist pure states $\varphi_1',\varphi_2'$, and
$\psi$ of $\UHF_d$, an increasing sequence $\{k_n\}$ in $\N$ and an increasing
sequence $\{B_n\}$ of finite type I subfactors of $A$ such that
\begin{eqnarray*}
&&\|\varphi_i-\varphi_i'\|<\varepsilon \\
&&\varphi_1'|B_{2n+1}\quad \mbox{is pure} \\
&&\varphi_2'|B_{2n}\quad \mbox{is pure} \\
&&\psi|B_{6k-1}\cap B_{6k-3}'=\varphi_1'|B_{6k-1}\cap B_{6k-3}' \\
&&\psi|B_{6k+2}\cap B_{6k}'=\varphi_2'|B_{6k+2}\cap B_{6k}' \\
&&\psi|B_{6k}\cap B_{6k-1}'\quad \mbox{is pure}, \\
&&\psi|B_{6k-3}\cap B_{6k-4}'\quad \mbox{is pure}, \\
&&k_{n+1}-k_n\to\infty \\
&&A_{k_1}\subset B_1\subset A_{k_2}\subset B_2\subset A_{k_3}\subset
        B_3\subset\cdots
\end{eqnarray*}
\end{lemma}

\begin{proof}
It follows from the previous lemma that there exist pure states $\varphi_i'$,
increasing sequences $\{B_{in}\}$ of finite type I subfactors of $A$, and an
increasing sequence $\{k_n\}$ in $\N$ such that
\begin{eqnarray*}
&&\|\varphi_i-\varphi_i'\|<\varepsilon\;, \\
&&\varphi_i|B_{in}\quad \mbox{is pure for}\ i=1,2\;, \\
&&A_{k_1}\subset B_{i1}\subset A_{k_2}\subset B_{i2}\subset
       A_{k_3}\subset\cdots
\end{eqnarray*}
By passing to subsequences of $\{k_n\}$ and $\{B_{in}\}$ and setting
$B_n=B_{1n}$ if $n$ is odd and $B_n=B_{2n}$ if $n$ is even, we may assume that
\begin{eqnarray*}
&&\varphi_1'|B_{2n+1}\quad \mbox{is pure} \\
&&\varphi_2'|B_{2n}\quad \mbox{is pure} \\
&&k_{n+1}-k_n\to\infty \\
&&A_{k_1}\subset B_1\subset A_{k_2}\subset B_2\subset A_{k_3}\subset\cdots \\
\end{eqnarray*}
Then $\varphi_1'$ has a tensor product decomposition into pure states on the
matrix subalgebras $B_{2n+1}\cap B_{2n-1}'$, and $\varphi_2'$ likewise on the
subalgebras $B_{2n}\cap B_{2n-2}'$. Thus we can define a pure state $\psi$ by
requiring that it decomposes under the tensor product decomposition
\begin{eqnarray*}
\lefteqn{\hspace*{-2em}\ldots\otimes(B_{6k-4}\cap B_{6k-6}')\otimes
     (B_{6k-3}\cap B_{6k-4}')\otimes(B_{6k-1}\cap B_{6k-3}')} \\
&&\otimes(B_{6k}\cap B_{6k-1}')\otimes(B_{6k+2}\cap B_{6k}')\otimes\cdots
\end{eqnarray*}
into states given by:
\begin{eqnarray*}
&&\psi|B_{6k-1}\cap B_{6k-3}'=\varphi_1'|B_{6k-1}\cap B_{6k-3}'\;, \\
&&\psi|B_{6k+2}\cap B_{6k}'=\varphi_2'|B_{6k+2}\cap B_{6k}'\;, \\
&&\psi|B_{6k}\cap B_{6k-1}'\quad \mbox{is an arbitrary pure state,} \\
&&\psi|B_{6k-3}\cap B_{6k-4}'\quad \mbox{is an arbitrary pure state}.
\end{eqnarray*}
\end{proof}

Recall that $\tau$ is the gauge action of $\T$ on ${\mathcal O}_d$, i.e.,
$$
\tau_z(s_i)=zs_i\;,\qquad z\in\T\;.
$$
Let $\varepsilon$ be the conditional expectation of ${\mathcal O}_d$ onto
$\UHF_d$ defined by
$$
\varepsilon(x)=\int_\T \tau_z(x)\frac{|dz|}{2\pi}\;,\qquad
x\in{\mathcal O}_d\;.
$$
Note that if $\varphi$ is a gauge-invariant state of ${\mathcal O}_d$, then
$$
\varphi=\varphi|_{\UHF_d}\circ\varepsilon\;.
$$
Recall that $\lambda$ is canonical endomorphism of ${\mathcal O}_d$:
$\lambda(x)=\sum\limits_{i=1}^d s_ixs_i^\ast$, $x\in{\mathcal O}_d$, and
that the restriction of $\lambda$ to $\UHF_d$ is the one-sided shift $\sigma$.

\begin{lemma}
If $\varphi$ is a gauge-invariant state on ${\mathcal O}_d$ then the
following conditions are equivalent:

{\rm (i)}\quad $\varphi$ is pure

{\rm (ii)} \
$\varphi|_{\UHF_d}$ is pure and

\qquad $\varphi|_{\UHF_d}\circ\sigma^n$ is disjoint from $\varphi$ for
$\;n=1,2,\ldots$
\end{lemma}

\begin{proof}
(i)$\Rightarrow$(ii). Since $\varphi$ is pure, and gauge-invariant, it
follows that $\varphi|_{\UHF_d}$ is pure. Let $p$ be the support projection
of $\varphi$ in ${\mathcal O}_d^{\ast\ast}$. Since $p$ is minimal, and
$\varphi$ is gauge-invariant, it follows that for any $a\in\UHF_d$ and
any multi-index $I=(i_1,i_2,\ldots,i_n)$ with $|I|=n\geq 1$,
$$
pas_I p=\varphi(as_I)p=0\;,
$$
where $s_I=s_{i_1}s_{i_2}\ldots s_{i_n}$. Thus we obtain that
$$
p(\UHF_d)\lambda^n(p)=0\;,
$$
which implies that $\varphi|_{\UHF_d}\circ \sigma^n$ is disjoint from
$\varphi$.
\bigskip

(ii)$\Rightarrow$(i). Let $p$ be the support projection of
$\varphi|_{\UHF_d}$ in
$\UHF_d^{\ast\ast}\subset{\mathcal O}_d^{\ast\ast}$. It suffices to show that
for any multi-indices $I, J$
$$
ps_Is_J^\ast p\in\C p
$$
since the linear span of $s_Is_J^\ast$ is dense in ${\mathcal O}_d$. If
$|I|\not=[J|$, we have that $ps_Is_J^\ast p=0$ by using the fact that
$\varphi|_{\UHF_d}\circ\sigma^n$ is disjoint from $\varphi$ for
$n=\Big||I|-|J|\Big|$. If $[I|=|J|$, we have that $ps_Is_J^\ast
p=\varphi(s_Is_J^\ast)p$ since $\varphi|_{\UHF_d}$ is pure.
\end{proof}

\begin{lemma}
Let $\varphi_1$ and $\varphi_2$ be gauge-invariant pure states of ${\mathcal
O}_d$ such that all
\newline
$\varphi_i|_{\UHF_d}\circ\sigma^n$, $i=1,2$,
$n=0,1,2,\ldots$ are mutually disjoint. Then there exists an automorphism
$\alpha$ of ${\mathcal O}_d$ such that
$\alpha\circ\tau_z=\tau_z\circ\alpha$, $z\in\T$ and
$\varphi_1=\varphi_2\circ\alpha$.
\end{lemma}

\begin{proof}
By Lemma~4, $\psi_1=\varphi_1|_{\UHF_d}$ and $\psi_2=\varphi_2|_{\UHF_d}$ are
pure states on $\UHF_d$. Applying Lemma~3 on $\psi_1,\psi_2$ in lieu of
$\varphi_1,\varphi_2$, with $\varepsilon=1$, we obtain pure states
$\psi_1',\psi_2'$ and $\psi$ of $\UHF_d$ with the properties given there.
Since $\psi_i$ is equivalent to $\psi_i'$,
$\varphi_i'=\psi_i'\circ\varepsilon$ is a pure state of ${\mathcal O}_d$ by
Lemma~4 and this state is equivalent to $\varphi_i=\psi_i\circ\varepsilon$.
By Kadison's transitivity theorem we have a unitary $u\in\UHF_d$ such that
$\psi_i'=\psi_i\circ\Ad u$; it follows that $\varphi_i'=\varphi_i\circ\Ad u$.

It is not automatical that $\psi$ satisfies the condition that all
$\psi\circ\sigma^n$, $n=0,1,2,\ldots$ are mutually disjoint and are disjoint
from $\psi_i'\circ\sigma^n$. But using the freedom in constructing
$\psi|_{B_{6k}\cap B_{6k-1}'}$ and $\psi|_{B_{6k-3}\cap B_{6k-4}'}$
successively, we can certainly impose this condition.

Thus we obtain three pure states $\psi_1',\psi_2',\psi$ of $\UHF_d$ such that
all $\psi_i'\circ\sigma^n$, $\psi\circ\sigma^n$ are mutually disjoint and
$\psi_i'$ and
$\psi$ are spotwise asymptotically equal as specified in Lemma~3. It now
suffices to prove the lemma for the pairs
$(\psi_1'\circ\varepsilon,\psi\circ\varepsilon)$ and
$(\psi_2'\circ\varepsilon, \psi\circ\varepsilon)$. Thus replacing
$\varphi_1,\varphi_2$ by one of these pairs, we may assume the lemma satisfy
the additional condition that there exists an increasing sequence $\{k_n\}$ in
$\N$ and an increasing sequence $\{B_n\}$ of finite type I subfactors of
$\UHF_d$ such that
\begin{eqnarray*}
&&A_{k_1}\subset B_1\subset A_{k_2}\subset B_2\subset A_{k_3}
       \subset B_3\subset \\
&&\varphi_i|_{B_{3n+1}}\quad \mbox{is pure}\;, \\
&&\varphi_1|_{B_{3n+3}\cap B_{3n+1}'}=\varphi_2|_{B_{3n+3}\cap
        B_{3n+1}'}\quad \mbox{is pure} \\
&&k_{3n+3}-k_{3n+2}\to\infty\;.
\end{eqnarray*}

We shall construct a sequence $\{v_n\}$ of unitaries in $\UHF_d$ such that
\newline
$\alpha=\lim\limits_{n\to\infty}\Ad(v_n v_{n-1}\ldots v_1)$ defines an
automorphism of ${\mathcal O}_d$ with $\varphi_1=\varphi_2\circ\alpha$. To
ensure the existence of the limit we choose the unitaries such that they
mutually commute and $\sum\|\lambda(v_n)-v_n\|<\infty$. Since $\alpha$
commutes with the gauge action $\tau$, this will complete the proof.

We fix a large $N\in\N$. We choose $n_1$ so large that the support
projections $e_i^{(1)}=\supp(\varphi_i|_{B_{3n_1+1}})$ are almost orthogonal
and
$k_{3n_1+3}-k_{3n_1+2}>2^{2(N+1)}$. Let $w_1$ be a partial isometry in
$B_{3n_1+1}$ with $w_1^\ast w_1=e_1^{(1)}$, $w_1w_1^\ast=e_2^{(1)}$. By the
polar decomposition of the approximate unitary
$$
w_1+(1-e_2^{(1)})w_1^\ast(1-e_1^{(1)})+(1-e_2^{(1)})(1-e_1^{(1)})\;,
$$
we obtain a unitary $v_1\in B_{3n_1+1}$ such that
$$
v_1 e_1^{(1)}=w_1
e_1^{(1)}=e_2^{(1)}w_1=e_2^{(1)}v_1\in B_{3n_1+1}
$$
and
$v_1(1-e_2^{(1)})(1-e_1^{(1)})\approx (1-e_2^{(1)})(1-e_1^{(1)})$.

We next choose $n_2>n_1$ so large that
$$
\sigma^n\circ\supp(\varphi_i|_{B_{3n_2+1}\cap B_{3n_1+3}'}),\quad
i\!=\!1,2,\;\,
n\!=\!-2^{N-1},-2^{-N+1}+1,\ldots,0,\ldots,2^{N+1}
$$
are almost orthogonal and
$k_{3n_2+2}-k_{3n_1+1}>2^{2(N+2)}$. (Though $\sigma$ is an endomorphism,
$\sigma^{-n}$ on $B_{3n_2+1}\cap B_{3n_1+3}'$ is well defined for
$n=1,2,\ldots,k_{3n_1+2}$.) Let $w_2$ be a partial isometry in
$B_{3n_2+1}\cap B_{3n_2+3}'$ such that
$$
w_2^\ast w_2=e_1^{(2)}=\supp(\varphi_1|_{B_{3n_2+1}\cap B_{3n_1+3}'})
$$
and
$$
w_2 w_2^\ast=e_2^{(2)}=\supp(\varphi_2|_{B_{3n_2+1}\cap B_{3n_1+3}'})\;,
$$
and let $\zeta$ be a partial isometry in $A_{k_{3n_2+2}+1}\cap
A_{k_{3n_1+3}}'$ such that $\zeta^\ast\zeta=e_1^{(2)}$ and
$\zeta\zeta^\ast=\sigma(e_1^{(2)})$.

Assume for the moment that $\sigma^\ell(e_i^{(2)})$, $i=1,2$;
$\ell=-2^{N+1},-2^{N+1}+1,\ldots,2^{N+1}$ are all orthogonal and set
$$
e_{ij}=
\left\{\!\!
\begin{array}{ll}
\sigma^{i-1}(\zeta)\sigma^{i-2}(\zeta)\ldots \sigma^{j}(\zeta)
        & \quad i>j \\ [.5ex]
\sigma^i(e_1^{(2)}) & \quad i=j \\ [.5ex]
\sigma^i(\zeta^\ast)\sigma^{i+1}(\zeta^\ast)\ldots
         \sigma^{j-1}(\zeta^\ast) & \quad i<j\end{array}\right.
$$
for $i,j=-2^{-N+1},\ldots,2^{N+1}$. Then $(e_{ij})$ is a family of matrix
units such that $\sigma(e_{ij})=e_{i+1,j+1}$ when $|i|, |i+1|,|j|, |j+1|\leq
2^{N+1}$. Let
\begin{eqnarray*}
\lefteqn{E=e_1^{(2)}+
\sum_{\ell=1}^{2^{N+1}-1}(1-e_1^{(2)})
\bigg\{ \frac{2^{N+1}-\ell}{2^{N+1}}e_{\ell,\ell}
     +\frac{\ell}{2^{N+1}}e_{\ell-2^{N+1},\ell-2^{N+1}}} \\
&&+\frac{1}{2^{N+1}}
     \sqrt{(2^{N+1}-\ell)\ell}\;
(e_{\ell,\ell-2^{-N+1}}+e_{\ell-2^{-N+1},\ell})\bigg\} (1-e_1^{(2)})
\end{eqnarray*}
as in \cite{Kis95}. Then $E$ is a projection in $D_2=A_{(k_{3n_2+2}+2^{N+1})}
\cap A_{(k_{3n_1+3}-2^{N+1})}'$ and satisfies
$$
\|\sigma(E)-E\|\sim \frac{1}{2^{\frac{N+1}{2}}}\;.
$$
Let $w=w_2+(1-e_2^{(2)})
\Big( \sum\limits_{\ell=1}^{2^{N+1}}(\sigma^\ell(w_2)+
\sigma^{-\ell}(w_2))\Big)(1-e_1^{(2)})$ and
$$
v=wE+(1-F)w^\ast(1-E)+(1-F)(1-E)
$$
where $F=wEw^\ast$.

By the orthogonality assumption on $\sigma^{\ell}(e_i^{(2)})$,
 $v$ is a unitary in
$D_2$ and satisfies
\begin{eqnarray*}
&&\|\sigma(v)-v\|\approx\|\sigma(E)-E\|\;, \\
&&v e_1^{(2)}=w_2e_1^{(2)}=e_2^{(2)}w_2=e_2^{(2)}v\;.
\end{eqnarray*}
Note also that $v$ commutes with $v_1$ and $e_i^{(1)}$.

Now, the projections $\sigma^{\ell}(e_i^{(2)})$, $i=1,2$,
$\ell=-2^{N+1},\ldots,2^{N+1}$ are not actually orthogonal but
choosing $n_2$ so large that they are very close to being orthogonal,
 we may obtain a unitary $v_2$ in $D_2$ by polar decomposition of $v$
such that $v_2$ satisfies the same conditions as above, i.e.,
\begin{eqnarray*}
&&v_2e_1^{(2)}=w_2e_1^{(2)}=e_2^{(2)}w_2=e_2^{(2)}v_2\in B_{3n_2+1}\cap
     B_{3n_1+3}'\;, \\
&&\|\lambda(v_2)-v_2\|\sim 2\raise.5ex\hbox{$^{-\frac{N+1}{2}}$}
\end{eqnarray*}
and $v_2\in D_2$.

Since
\begin{eqnarray*}
\lefteqn{
\supp(\varphi_1|_{B_{3n_2+1}})} \\
&&=\supp(\varphi_1|_{B_{3n_1+1}})
    \supp(\varphi_1|_{B_{3n_1+3}\cap B_{3n_1+1}'})
    \supp(\varphi_1|_{B_{3n_2+1}\cap B_{3n_1+3}'}) \\
&&=
e_1^{(1)}pe_1^{(2)}
\end{eqnarray*}
with $p=\supp(\varphi_1|_{B_{3n_1+3}\cap B_{3n_1+1}'})=
\supp(\varphi_2|_{B_{3n_1+3}\cap B_{3n_1+1}'})$, and since the operators
$v_1e_1^{(1)}=e_2^{(1)}v_1,p$, and $v_2e_1^{(2)}=e_2^{(2)}v_2$ commute, we
obtain that
\begin{eqnarray*}
\lefteqn{v_1v_2\cdot\supp(\varphi_1|_{B_{3n_2+1}})=v_1v_2e_1^{(1)}
        pe_1^{(2)}} \\
&&=v_1e_1^{(1)}v_2e_1^{(2)}p \\
&&=e_2^{(1)}v_1e_2^{(2)}v_2p \\
&&=pe_2^{(1)}e_2^{(2)}v_1v_2=\supp(\varphi_2|_{B_{3n_2+1}})v_1v_2\;.
\end{eqnarray*}
Here we have also used the fact that $v_1$ commutes with $e_2^{(2)}$.
We repeat this procedure. Thus we obtain an increasing sequence $\{n_k\}$ in
$\N$ and a sequence $\{v_k\}$ of mutually commuting unitaries such that
\begin{eqnarray*}
&&\|\lambda(v_k)-v_k\|\sim
     2\raise.5ex\hbox{$^{-\frac{N+k}{2}}$}\;, \\
&&v_k e_1^{(k)}=e_2^{(k)}v_k\in{\mathcal B}_{3n_k+1}\cap
{\mathcal B}_{3n_{k-1}+3}'
\end{eqnarray*}
where
$$
e_i^{(k)}=\supp(\varphi_i|_{{\mathcal B}_{3n_k+1}\cap
{\mathcal B}_{3n_{k-1}+3}'})\;,
$$
and such that $\Ad(v_k\ldots v_1)$ maps
$\supp(\varphi_1|_{{\mathcal B}_{3n_k+1}})$ into
$\supp(\varphi_2|_{{\mathcal B}_{3n_k+1}})$. Then the limit
$\alpha=\lim\limits_k\Ad(v_k\ldots v_1)$ defines the desired automorphism.
\end{proof}

\begin{theorem}
Let $\varphi_1$ and $\varphi_2$ be gauge-invariant pure states of ${\mathcal
O}_d$. Then there exists an automorphism $\alpha$ of ${\mathcal O}_d$ such
that $\varphi_1=\varphi_2\circ\alpha$.
\end{theorem}

\begin{proof}
If $\varphi_1$ is disjoint from $\varphi_2$, then it follows that
$(\varphi_i|_{\UHF_d})\circ \sigma^n=\varphi_i\circ\lambda^n|_{\UHF_d}$,
$i=1,2$,
$n=0,1,2,\ldots$ are mutually disjoint (by Lemma~4); thus the assertion
follows from Lemma~5. If
$\varphi_1$ is equivalent to
$\varphi_2$, there is a unitary $u\in{\mathcal O}_d$ such that
$\varphi_1=\varphi_2\Ad u$ (by Kadison's transitivity).
\end{proof}

\section{Pure states mapped into Cuntz states by endomorphisms}

\setcounter{equation}{0}

There is a one-to-one correspondence between the set ${\mathcal U}({\mathcal
O}_d)$ of unitaries of ${\mathcal O}_d$ and the set $\End({\mathcal O}_d)$ of
unital endomorphisms of ${\mathcal O}_d$; if $u\in {\mathcal U}({\mathcal
O}_d)$, the endomorphism $\alpha_u$ is defined by $\alpha_u(s_i)=us_i$ and if
$\alpha\in\End({\mathcal O}_d)$, $\alpha$ corresponds to the unitary $u$
defined by $u=\sum\limits_{i=1}^d \alpha(s_i)s_i^\ast$. Define
\begin{eqnarray*}
&&{\mathcal U}_i=\{u\in {\mathcal U}({\mathcal O}_d)|\quad \alpha_u\;\,
      \mbox{is an inner automorphism}\} \\
&&{\mathcal U}_a=\{u\in {\mathcal U}({\mathcal O}_d)|\quad \alpha_u\;\,
     \mbox{is an automorphism}\} \\
&&{\mathcal U}_s={\mathcal U}({\mathcal O}_d)\setminus {\mathcal U}_a\;.
\end{eqnarray*}

\begin{proposition}
Let ${\mathcal U}_i,{\mathcal U}_a,{\mathcal U}_s$ be as above.

(i)\quad
${\mathcal U}_i$ is a dense subset of ${\mathcal U}({\mathcal O}_d)$.

(ii)\quad
${\mathcal U}_a$ is a dense $G_\delta$ subset of ${\mathcal U}({\mathcal
O}_d)$.

(iii)
$\;{\mathcal U}_s$ is a dense $F_\sigma$ subset of ${\mathcal U}({\mathcal
O}_d)$.
\end{proposition}

\begin{proof}
M.~R\o{}rdam proved (i) in [R\o{}r93] and the other statements are more or
less known.

We shall give a proof of (i). We again denote by $\lambda$ the canonical
endomorphism of ${\mathcal O}_d:\lambda(x)=\sum\limits_{i=1}^d s_ixs_i^\ast$,
$x\in{\mathcal O}_d$. Since the unitary corresponding to $\Ad v$ is
$v\lambda(v^\ast)$, it suffices to show that $v\lambda(v^\ast)$, $v\in
{\mathcal U}({\mathcal O}_d)$, is dense in ${\mathcal U}({\mathcal O}_d)$. If
$\UHF_d$ denotes the C*-subalgebra generated by $s_{i_1}s_{i_2}\ldots
s_{i_n}s_{j_n}^\ast\ldots s_{j_1}^\ast$, then we mentioned in the
introduction that $\UHF_d$ is isomorphic to the UHF algebra
$\bigotimes\limits_{\N}M_d$ and $\lambda|\UHF_d$ corresponds to the one-sided
shift on $\bigotimes\limits_{\N}M_d$. Thus $\lambda|\UHF_d$ satisfies the
Rohlin property, \cite{BKRS93}, \cite{Kis95}. In particular for any $n$ and
$\varepsilon>0$ there is an orthogonal family $e_0,e_1,\ldots,e_{n-1}$ of
projections in
$\UHF_d$ such that
\begin{eqnarray*}
&&\sum_{i=0}^{d^n-1} e_i=1 \\
&&\|\lambda(e_i)-e_{i+1}\|<\varepsilon
\end{eqnarray*}
with $e_{d^n}=e_0$. The similar properties hold for $\Ad u\circ\lambda$,
i.e., if $\UHF_d^u$ denotes the C*-subalgebra generated by
$us_{i_1}us_{i_2}\ldots us_{i_n}s_{j_n}^\ast u^\ast\ldots
s_{j_1}^\ast u^\ast$, then $\Ad u\circ\lambda|\UHF_d^u$ corresponds to the
one-sided shift on $\bigotimes\limits_{\N}M_d$. Hence for any $n$ and
$\varepsilon>0$ there is an orthogonal family $f_0,f_1,\ldots,f_{d^n-1}$ of
projections in $\UHF_d^u$ such that
\begin{eqnarray*}
&&\sum_{i=0}^{d^n-1}f_i=1 \\
&&\|\Ad u\circ\lambda(f_i)-f_{i+1}\|<\varepsilon
\end{eqnarray*}
with $f_{d^n}=f_0$. Suppose we have chosen such projections $e_i,f_i$ for
the same $n$. Since $K_0({\mathcal O}_d)=\Z/(d-1)\Z$, we have that
$[e_0]=1=[f_0]$ in $K_0({\mathcal O}_d)$ and so obtain a partial isometry
$w\in{\mathcal O}_d$ such that $w^\ast w=e_0$, $ww^\ast=f_0$. We find
unitaries $v_1,v_2\in{\mathcal O}_d$ such that $\Ad
v_1\lambda(e_i)=e_{i+1}$, $\Ad v_2\Ad u\lambda(f_i)=f_{i+1}$, and
$\|v_1-1\|\approx 0$, $\|v_2-1\|\approx 0$ (depending on $\varepsilon$). Let
$$
z=w^\ast(L_{v_2u}R_{v_1^\ast}\lambda)^{d^n}(w)
$$
where $R_{v_1^\ast}$ is the right multiplication by $v_1^\ast$ and
$L_{v_2u}$ is the left multiplication by $v_2u$. Since
$(L_{v_2u}R_{v_i^\ast}\lambda)^i(w)$ is a partial isometry with initial
projection $e_i$ and final projection $f_i$, $z$ is a unitary in
$e_0{\mathcal O}_d e_0$. Since $K_1({\mathcal O}_d)=0$ and ${\mathcal O}_d$
has real rank zero, we find a sequence
$z_0,z_1,\ldots,z_{d^n-1}$ of unitaries in $e_0{\mathcal O}_d e_0$ such that
$z_0=z$, $z_{d^n-1}=1$,
$$
\|z_i-z_{i+1}\|<4/d^n\;.
$$
Define a unitary $v$ by
$$
v=\sum_{i=0}^{d^n-1}(L_{v_2u}R_{v_1^\ast}\lambda)^i(wz_i)
$$
Then since
\begin{eqnarray*}
\lefteqn{v-(L_{v_2u}R_{v_1^\ast}\lambda)(v)} \\
&&=\sum_{i=1}^{d^n-1}(L_{v_2u}R_{v_1^\ast}\lambda)^i(wz_i-wz_{i-1})
        +wz_0-(L_{v_2u}R_{v_1^\ast}\lambda)^{d^n}(w)\;,
\end{eqnarray*}
it follows that
$$
\|v-L_{v_2u}R_{v_1^\ast}\lambda(v)\|<4/d^n
$$
or
$$
\|v-u\lambda(v)\|\lesssim 4/d^n\;.
$$
This completes the proof of (i).

Since ${\mathcal U}_a\supset {\mathcal U}_i$, ${\mathcal U}_a$ is dense. That
${\mathcal U}_a$ is a
$G_\delta$ set follows from
$$
{\mathcal U}_a=\bigcap_n \bigcap_j \bigcup_i \bigg\{
u\in {\mathcal U}({\mathcal O}_d);
\|\alpha_u(x_i)-x_j\|<\frac{1}{n}\bigg\}
$$
where $\{x_i\}$ is a dense sequence in ${\mathcal O}_d$.

If ${\mathcal U}_a$ contains a non-empty open set, then it follows that
${\mathcal U}_a={\mathcal U}({\mathcal O}_d)$ or ${\mathcal U}_s=\emptyset$.
Because for any unitaries $u,w$ of ${\mathcal O}_d$ we find a unitary $v$
such that
$w\lambda(v)\approx vu$. (Apply the previous argument for the endomorphism
$\Ad u\circ\lambda$ instead of
$\lambda$ and the unitary $wu^\ast$.) Since
$v{\mathcal U}_a\lambda(v^\ast)={\mathcal U}_a$ for any unitary
$v\in{\mathcal O}_d$, the above fact implies that ${\mathcal U}_a$ contains an
arbitrary unitary. But we know that
${\mathcal U}_s\not=\emptyset$. For example if
$u=\sum s_is_js_i^\ast s_j^\ast$, then $\alpha_u=\lambda$ and
$\lambda({\mathcal O}_d)'\simeq M_d$. Thus we obtain that ${\mathcal U}_s$ is
dense.
\end{proof}

For a unit vector $\xi\in\C^d$ we have defined the Cuntz state $f_\xi$ of
${\mathcal O}_d$ by
$$
f_\xi(s_{i_1}\ldots s_{i_m}s_{j_n}^\ast\ldots
s_{j_1}^\ast)=\xi_{i_1}\ldots
\xi_{i_m} \overline{\xi_{j_n}}\ldots\overline{\xi_{j_1}}
$$
It follows that $f_\xi$ is a unique pure state of ${\mathcal O}_d$ satisfying
$$
f_\xi\bigg(\sum_{i=1}^d\overline{\xi_i}s_i\bigg)=1\;.
$$
Let $F$ be the linear span of $s_is_j^\ast$, $i,j=1,\ldots,d$. Then $F$ is
isomorphic to $M_d$ and each unitary $u$ in $F$ defines an automorphism
$\alpha_u$ of ${\mathcal O}_d$. This group of automorphisms acts transitively
on the compact set of Cuntz states.

We denote by $f_0$ the Cuntz state $f_\xi$ with $\xi=(1,0,\ldots,0)$.

\begin{proposition}
If $\varphi$ is a pure state of ${\mathcal O}_d$, there is a unital
endomorphism $\alpha$ of ${\mathcal O}_d$ such that $\varphi\circ\alpha=f_0$,
where $f_0$ is the Cuntz state defined above. Furthermore $\alpha$ may be
chosen so that $\pi_\varphi\circ\alpha({\mathcal O}_d)''$ contains the
one-dimensional projection onto $\C\Omega_\varphi$.
\end{proposition}

\begin{proof}
It suffices to show that if $\varphi$ is a pure state there is a unitary
$u\in{\mathcal O}_d$ such that
$$
\varphi(us_1)=1\;.
$$
Since ${\mathcal O}_d$ has real rank zero, there is a decreasing sequence
$(e_n)$ of projections in ${\mathcal O}_d$ such that $\varphi$ is the unique
state satisfying $\varphi(e_n)=1$ for $n=1,2,\ldots$, i.e., $(e_n)$ converges
to the support projection of $\varphi$ in ${\mathcal O}_d^{\ast\ast}$. We may
further assume that $[e_n]=0$ in $K_0({\mathcal O}_d)$.

Pick up a projection $e=e_n$ such that $\varphi(e)=1$ and $e<1$. Then
$es_1^\ast$ is a partial isometry with initial projection $s_1es_1^\ast$ and
final projection $e$. Let $w$ be a partial isometry such that $w^\ast
w=1-s_1es_1^\ast$ and $ww^\ast=1-e$. Then
$u=es_1^\ast+w$ is a unitary in ${\mathcal O}_d$ such that
$$
us_1e=(es_1^\ast+w)s_1e=e\;.
$$
Thus we have that $\varphi(us_1)=1$.

To prove the last statement we shall modify $u$ so that $\varphi$ is the
unique state satisfying
$$
\varphi(us_1)=1\;.
$$
We have chosen $e=e_n$. We let
$$
h=\sum_{k=1}^\infty 2^{-k}e_{n+k}\;.
$$
Then $h$ is self-adjoint with $0\leq h\leq 1$ and $\varphi$ is the only state
satisfying $\varphi(h)=1$. Let
$$
u_1=e^{2\pi ih}u\;.
$$
Then $u_1 s_1e=e^{2\pi ih}e$ and the assertion follows.
\end{proof}

\end{document}